\documentclass[journal]{IEEEtran}
\usepackage{amssymb,times}
\usepackage{cite, epsfig}
\usepackage{ verbatim, amsfonts,graphicx,enumerate}
\usepackage{bm,array}
\usepackage[cmex10]{amsmath}
\usepackage{bm}
\usepackage{framed}
\usepackage{graphicx}
\usepackage{enumitem}
\usepackage{epstopdf}
\usepackage{amsthm} 
\usepackage{thmtools, thm-restate}

\usepackage{caption}
\usepackage{subcaption}
\usepackage{color} 
\usepackage{makecell}
\usepackage{enumitem}
\usepackage{algorithm}
\usepackage{varwidth}
\usepackage{algpseudocode}
\usepackage{tikz}
\theoremstyle{plain}
\newtheorem{thm}{Theorem}

\theoremstyle{definition}

\theoremstyle{remark}

\def\h{\mathbf{h}}
\def\v{\mathbf{v}}

\def\x{\mathbf{x}}
\def\s {{\mathbf{s}}}
\def\d {{\mathbf{d}}}

\def\w{\mathbf{w}}
\def\y{\mathbf{y}}

\def\E{\mathbf{E}}

\def \A {\mathbf{A}}

\def \S {\mathcal{S}}

\def \bell	{{\boldsymbol{\ell}}}

\providecommand{\norm}[1]{\left\|#1\right\|}

\providecommand{\norm}[1]{\left\|#1\right\|}
\title{An Online Approach to D2D Trajectory Utility Maximization Problem}
	\author{\IEEEauthorblockN{Amrit S. Bedi,$^{\star}$ Ketan Rajawat,$^{\star}$ and  Marceau Coupechoux$^{\dagger}$}\\
	\IEEEauthorblockA{$^{\star}$Department of Electrical Engineering, Indian Institute of Technology Kanpur, India 208016\\}
	\IEEEauthorblockA{$^\dagger$LTCI, Telecom ParisTech, University Paris-Saclay 75013, Paris, France \\E-mail:$^{\star}$\{amritbd, ketan\}@iitk.ac.in, $^{\dagger}$marceau.coupechoux@telecom-paristech.fr} \thanks{{Work in this paper was supported by the Indo-French Centre for the Promotion of Advanced Research (IFCPAR/CEFIPRA).} }}

\begin{document}
		\maketitle
\begin{abstract}	
This paper considers the problem of designing the user trajectory in a device-to-device communications setting. We consider a pair of pedestrians connected through a D2D link. The pedestrians seek to reach their respective destinations, while using the D2D link for data exchange applications such as file transfer, video calling, and online gaming. In order to enable better D2D connectivity, the pedestrians are willing to deviate from their respective shortest paths, at the cost of reaching their destinations slightly late. A generic trajectory optimization problem is formulated and solved for the case when full information about the problem in known in advance.

Motivated by the D2D user's need to keep their destinations private, we also formulate a regularized variant of the problem that can be used to develop a fully online algorithm. The proposed online algorithm is quite efficient, and is shown to achieve a sublinear \emph{offline} regret while satisfying the required mobility constraints exactly. The theoretical results are backed by detailed numerical tests that establish the efficacy of the proposed algorithms under various settings. 
	
\end{abstract}
\section{Introduction}
The emergence of social media and the associated mobile applications has ushered a culture of constant connectivity.  With the advent of data intensive applications, the service providers now face the challenge of providing seamless connectivity `on the go' \cite{li2009future,chin2014emerging}. The need to maintain high spectral efficiency has prompted the researchers to look beyond the traditional cellular architecture and develop innovative features such as device-to-device (D2D) communications \cite{lin2000multihop}. Enabling direct communications between nearby cellular users has the potential to not only improve spectrum utilization, throughput, and energy efficiency, but also enable disruptive peer-to-peer applications and services \cite{ji2016wireless,archi_D2D,yu2014joint}. Indeed, the D2D paradigm allows the users to sustain connectivity at low costs and without overloading the cellular network \cite{asadi2014survey}. 

Thanks to the content heavy nature of the modern social media platforms, the cybercitizens of today are increasingly willing to modify their behavior in order to stay connected. For instance, urban areas have seen an unprecedented increase in the number of wifi hotspots \cite{mota2013feasibility}. As more D2D devices surface, it is likely that the users will be willing to modify their daily commute so as to stay connected. For instance, pedestrians may be willing to increase their commute times by a fraction so as to not only exchange more data and connect longer, but also sustain a higher quality of service. 

This paper formulates the D2D trajectory optimization problem for commuting users, where the goal is to maximize a user-specific utility function depending on the download data rate or of the signal strength, while ensuring that the users reach their destinations. Specifically, we consider a pair of mobile users following a certain trajectory, say from a starting point to a destination. The users exchange data on a D2D link and need to sustain high data rates. To this end, the users are willing to take detour from the most direct path from the starting point to their respective destinations. The goal now is to design an optimal trajectory that the users must take so as to not only maximize their data rates but also reach their respective destinations within a preset delay. 

This paper considers the trajectory design problem from an optimization perspective. It is shown that with the current formulation, finding the optimal trajectory is relatively straightforward if both the users share their starting points and destinations in advance. However, such an offline approach is not only impractical but also potentially insensitive to the privacy requirements of a user. The need for the users to hide their destinations prompts us to consider the problem from the perspective of a single user. Within this privacy-preserving trajectory optimization framework, the two users still connect via the D2D link but do not reveal their future trajectories to each other. Consequently, each user designs its own trajectory in an online manner while taking into account the uncertainity in the trajectory of the other user. 

The problem of path planning or trajectory optimization has previously been considered within the context of rolling horizon planning, model predictive control (MPC), and online learning \cite{falcone2007linear,kim2014model,mokhtari2016online}. To this end, we propose an MPC-inspired algorithm that is shown to perform close to the offline problem but is prone to an infeasibility problem. Next, in order to circumvent this issue, we approach the problem from an online learning perspective. We advocate a modified version of the online gradient descent (OGD)  algorithm that can not only handle generic cost functions, but also incurs a sublinear dynamic regret compared to the offline problem. The resulting regret bound is also the strongest such result for generic online learning problems, significantly advancing the state-of-the-art. The subsequent application of the proposed algorithm to the trajectory optimization problem at hand provide insightful results and are verified via detailed simulations. 

In summary, the present works makes two major contributions, namely, (a) the D2D trajectory optimization problem is formulated and solved using novel online algorithms, and (b) a modified OGD algorithm is proposed that is shown to incur a sublinear offline regret and can be used to solve generic trajectory optimization problems. 
 
\subsection{Related work and contributions}
The problem of trajectory planning for a moving robot in a specified area has been widely studied in the robotics and control literature \cite{thomas2003towards,paden2016survey,goerzen2010survey,vasquez2016novel}. The motion planning problems are often categorized into two settings: \emph{static} and \emph{time varying} \cite{paden2016survey}.

For the time-varying setting that is of interest here, the standard approaches are based upon graph theory \cite{tsitsiklis1995efficient,takei2010practical,reif2000nonuniform,lavalle1998optimal,lavalle2006planning,likhachev2009planning,ferguson2006using}, Markov decision process based techniques \cite{fakoor2016humanoid,spaan2004point}, and model predictive control \cite{garcia1989model,falcone2007predictive,falcone2007linear,kim2014model}. Other advances in the path planning and robotics areas include dealing with uncertainties, dynamics, or multiple robots \cite{van2011lqg,liu2011coordinated,plaku2010motion}. The standard approach here is to generally formulate a complicated optimization problem and solve it using interior point method at each time step $t$. In constrast, the present paper specifically targets the problem from an online learning perspective, leading to a significantly simpler algorithm that is also amenable to well-defined performance guarantees.

In a similar vein, the path planning problem has also been considered under the aeges of reinforcement learning. Within this context, the problem is formulated as that of finding a sequence of feasible actions that take a robot from a source to destination \cite{al2013wind,konar2013deterministic}. Different from these works, { we do not assume that the distance to the next position from the current position is known to the moving user.}

The path planning problem has also been considered within the control theory literature \cite{chen2013retracted,barkaoui2014information}, where the goal is to design asymptotically stable systems. In contrast, the utility maximization framework considered is finite-horizon and consequently uses different tools for analysis and development of performance guarantees. 

The problem of trajectory design has also been considered from a variational perspective, with the trajectory given by a continuous time function. The variational problems are often solved numerically via discretization and consequently are not amenable to the regret bounds such as those developed here \cite{fraichard1998trajectory}.

Finally, the related problem of target tracking has been considered in the context of signal processing as well in online learning \cite{derenick2009optimal,bedi2017adversarial}. Within this context, it is common to assess the performance of online algorithms using the notion of dynamic regret, that quantifies the difference between the the cost achieved by the online algorithm and that achieved by an adaptive adversary \cite{mokhtari2016online}. The present work builds upon this formulation and utilizes a stronger notion of offline regret, wherein the adversary also has access to the information from the future \cite{chen2017online}. Different from the existing literature on dynamic regret, the present work provides a sublinear characterization of the offline regret \cite{chen2017online}.

\textbf{Notations:} All the scalars are represented by regular font and vectors by bold font. 

\section{Generic problem formulation}
This section formulates the D2D trajectory optimization problem. Consider a pair of users located on the $\mathbb{R}^2$ plane and connected via a D2D link. The location of the two users at discrete time $t \in \mathbb{N}$ is denoted by $\x_1(t), \x_2(t) \in \mathcal{S} \subset \mathbb{R}^2$, where $\mathcal{S}$ is the set of viable user locations. The origin and the destination of the $i$-th user are denoted by $\s_i$ and $\d_i$ respectively. For the base case, the $i$-th user does not modify its behavior and travels from $\s_i$ to $\d_i$ along the shortest path. When there are no obstacles, the shortest path is simply the straight line joining $\s_i$ and $\d_i$ and has length $\norm{\s_i-\d_i}_2$. On the other hand, when the user is only allowed to move along a grid, the length of the shortest path is the city-block or Manhattan distance $\norm{\s_i-\d_i}_1$. The distance between the source and the destination is henceforth denoted by $\norm{\s_i-\d_i}$, where the norm could be Euclidean, Manhattan, or any other convex distance metric. User $i$ travels at the maximum speed of $v_i$ units per time slot and therefore takes time $T_i:=\norm{\s_i-\d_i}/v_i$ to reach the destination via the shortest path. 

The two users communicate on the D2D link and derive a utility of $U(\norm{\x_1(t)-\x_2(t)}_2)$ at time $t$. Here, $U$ is a non-increasing function of $\norm{\x_1(t)-\x_2(t)}_2$. Examples of utility functions may include the average received signal strength modeled as 
\[
U_{\text{RSS}}(\norm{\x_1(t)-\x_2(t)}_2) = RSS =\frac{1}{\norm{\x_1(t)-\x_2(t)}_2^\alpha}
\]
where $\alpha$ is the path loss parameter, or functions thereof. For instance the average signal-to-noise ratio for the additive white gaussian noise channel with noise power $\sigma^2$ given by
\[
U_\text{SNR}(\norm{\x_1(t)-\x_2(t)}_2)= SNR = \frac{RSS}{\sigma^2+RSS}
\]
and the channel capacity for a channel with bandwidth $W$ given by 
\[
U_C(\norm{\x_1(t)-\x_2(t)}_2)=W\log_2(1+SNR)
\]
are some common examples. It is emphasized that in practice, the exact D2D rate would likely depend on a number of other factors, such as the application used, overheads, shadowing, etc. Given the inherent uncertainty in estimating the D2D rate in an urban setting, it may not necessarily be prudent to choose a complicated utility function. We set $T :=\min\{T_1,T_2\}$ and assume that the users disconnect as soon as one of them reaches the destination. Therefore, the cumulative utility obtained by the two users is given by $\sum_{t = 1}^{T} U(\norm{\x_1(t)-\x_2(t)}_2)$. 

We consider a scenario where the users are willing to take a longer path to their destinations in order to achieve a higher utility. Specifically, let $\delta_i$ be the excess delay user $i$ is willing to incur. In other words, it is acceptable for user $i$ to reach its destination in time $T_i+\delta_i$ so as to achieve a higher cumulative utility. Note that in the current formulation, the excess delay $\delta_i$ is an exogenous variable set by the user prior to starting. The D2D trajectory optimization problem can therefore be written as
\begin{subequations}\label{main}
	\begin{align}
	&&&\hspace{-5cm}\max_{\{\{\x_i(t)\}_{t=1}^{T_i+\delta_i}\}_{i=1}^2} \hspace{-5mm}\sum_{t = 1}^{\min_i\{T_i+\delta_i\}} U(\norm{\x_1(t)-\x_2(t)}_2) \nonumber\\
	\text{s. t. } \ \x_i(1)&=\s_i & i &=1, 2 \label{start_rel}\\
	\x_i(T_i+\delta_i)&=\d_i & i &= 1, 2  \label{dest_rel}\\
	\norm{\x_i(t+1)-\x_i(t)}&\leq v_i  & i &= 1,2\nonumber\\
	&&1\leq t &< T_i+\delta_i,\label{vel} \\
	\x_i(t) &\in \mathcal{S} & i  &=1, 2 \label{set}
	\end{align}
\end{subequations}
Here, the constraints in \eqref{start_rel} and \eqref{dest_rel} ensure that the user trajectories begin at their respective starting points and end at their corresponding destinations. The constraint in \eqref{vel} enforces the maximum velocity constraint under the appropriate norm. Finally, the constraint in \eqref{set} ensures that the user trajectories stay within the viable region. In an urban setting, the set $\S$ may represent the public areas such as roads, streets, alleys, etc. Such a constraint can generally be represented as a union of $K$ convex polygons, where each polygon represents an unobstructed area, such as the length and width of a road. Formally, each polygon is described by a set of affine inequalities of the form $\A_k\x_i(t) \leq \mathbf{b}_k$ for $1\leq k \leq K$ \cite{blackmore2006probabilistic}. It is emphasized that $\S$ is a union, not an intersection, of such polygons, and is generally non-convex.

The optimization problem in \eqref{main} is non-convex and consequently difficult to solve in general. While it may be possible to define an appropriate utility function that is concave in $\{\x_i(t)\}$, handling the non-convex set $\S$ is not straightforward and necessitates some simplifications. Indeed, a large subset of literature on path planning is dedicated to navigating obstacles through the use of heuristics. Within this context, a common approach towards handling \eqref{set} is discretization, wherein $\S$ is represented as a graph $G = (V,E)$. In particular, the user locations are restricted to the nodes $V$ of the graph, and the edges encode the movement cost between neighboring locations. The resulting problem is no longer a continuous domain optimization problem, but instead requires tools from graph theory. It can be observed that the complexity of the resulting problem now depends on the granularity of the discretization process, and can be high for moderate sized problems. 

This paper considers the continuous domain trajectory optimization problem that will subsequently be used to develop offline and online optimization algorithms. To this end, we consider a special case of \eqref{main} with the following two assumptions (a) the utility function $U$ is concave in $\x_1(t)$ and $\x_2(t)$ for all $t$; and (b) the set of viable locations $\S$ is convex. The first assumption is satisfied, for instance, if the function $U(\cdot)$ is concave and non-increasing, e.g., $U(x) = -x^2$. On the other hand, the second assumption is not generally satisfied in areas with obstacles, but is required in order to develop meaningful theoretical guarantees. As will be shown later, the proposed online algorithm may still be run with non-convex $\S$ if a projection operation can be calculated easily. However, developing performance guarantees for non-convex $\S$ is significantly harder and will not be pursued here. 

Having formulated the problem at hand, the subsequent section details the offline  approaches for solving \eqref{main}. 

\section{Offline approaches} \label{offline}
When $\S$ is convex and can be expressed as an intersection of a few half-spaces or through other convex functions, it is possible to solve \eqref{main} using an interior point algorithm. As discussed earlier, in order to solve \eqref{main}, it is necessary for the users to cooperate with each other and share the information about their final destinations prior to starting. The complexity of solving \eqref{main} depends on the granularity of time discretization, since the number of optimization variables is $T_1+T_2$.  

\begin{figure}
	\centering
	\hspace{0mm}\includegraphics[width=0.89\linewidth, height = 0.5\linewidth]{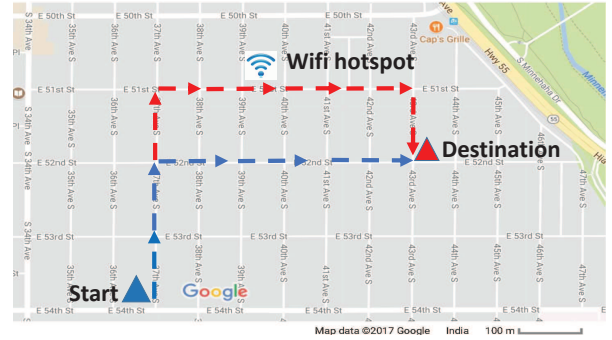}
	\caption{A simple fixed hotspot model, blue dashed line shows a direct path while red dashed line shows an alternative path to increase data rate.}
	\label{system_model}
	\vspace{-5mm}
\end{figure}

We remark of a special case where one of the users is simply a fixed hot-spot, i.e., $\x_2(t) = \h$. In this case, only the hotspot location is required in advance. Consequently, the user may calculate its desired route by simply selecting an appropriate value of the excess delay $\delta$. Such a use case is depicted in Fig.~\ref{system_model}.  It is evident that if $\delta$ is sufficiently large, it may be optimal to simply reach to the hot spot, wait as long as possible, and subsequently head to the destination at the maximum speed. 

In the general D2D setting however, the offline formulation in \eqref{main} has limited applicability. To begin with, the formulation encodes a static scenario, where both users cooperate and exchange their destinations prior to the starting. In practice, the destination of the user may change after the user has already started, rendering the solution to \eqref{main} useless. More importantly, the users may like to interact online but not necessarily be willing to reveal their final destinations due to privacy and security concerns. The a priori unavailability of the full problem information motivates the need for online algorithms that are capable of handling time-varying and uncertain parameters. 

To this end, consider the case when $\x_2(t)$ is an exogenous variable that is revealed at each time $t$. Likewise, we allow the final destination to be a time-varying quantity $\d_i(t)$. The resulting time-varying problem can no longer be expressed in the form of \eqref{main}. Such problems have traditionally been considered within the context of rolling horizon planning or model predictive control. The idea here is to solve the offline problem at each time $t$ repeatedly, every time a problem parameter changes. 

Here, we propose an MPC-like online algorithm (c.f. Algo.~\ref{algo1}) that can handle time-varying $\x_2(t)$ and $\d_2(t)$ in real time. In the absence of future information, we purpose to solve the full optimization problem from $\tau=1$ to $T+\delta$ at each time instant $t$. The idea here is to assume the current position of the exogenous user $\x_2(t)$ remains fixed for all $\tau=t, \ldots, T+\delta$ and plan accordingly. The step by step implementation of the proposed algorithm is described in Algo.~\ref{algo1}.  
\begin{algorithm}
	\caption{: MPC based algorithm }\label{algo1}
	\begin{algorithmic}[1]
		\item {\bf Given} $\s_1$, $\d_1$, $\nu_1$, $\hat{\x}_1=\s_1$, and $\x_2(1)$
		\item {\textbf{for} $t=1$ to $(T+\delta-1)$}
		\item \ \ \ \ \textbf{Solve} the following  optimization problem
		\ \ \ \small\begin{subequations}\label{main22}
			\begin{align}
			\{\x^M_1(\tau)\}_{\tau=1}^{T+\delta}:=\min_{\{\x_1(\tau)\}_\tau}  \sum_{\tau=1}^{T+\delta}&\left\| \x_1(\tau)-\x_2(t)\right\|^2 \\
			\text{s. t. }  \x_1(1)&=\s \label{cosnt12}\\
			\x_1(T+\delta)&=\d_1 \label{cosnt22}\\
			\left\| \x_1(\tau+1)-\x_1(\tau) \right\| &\leq v_1 \nonumber\\
			&\text{for }~\tau\in[1,T+\delta-1]\label{cosnt32}\\
			\x_1(\tau)&=\hat{\x}_1(\tau) \ \text{for} \ \tau=1\ \text{to}\ t.  \label{only_MPC}  				\end{align}
		\end{subequations}
		\item  \ \ \ \ \ Move to position $\hat{\x}_1(t+1)=\x_1^M(t+1)$. 
		\item \textbf{end} 
	\end{algorithmic}
\end{algorithm}

An important practical issue that arises for this class of algorithms is the possibility of the problem becoming infeasible at some time $t \geq 1$.  The problem could become infeasible, for instance, if $\d(t)$ changes abruptly at some point. Intuitively, the change in $\d(t)$ could be such that it is no longer possible for the user to reach the destination even when traveling at the maximum velocity towards $\d(t)$. The numerical tests reveal that for Algo.~\ref{algo1} to yield a feasible solution for all $1\leq t < T+\delta$, it is necessary that both $\norm{\d(t)-\d(t-1)}$ and $\norm{\x_2(t)-\x_2(t-1)}$ be small for all $t$. Moreover, analytically characterizing the problem parameters that yield feasible solutions in Algo.\ref{algo1} is not straightforward. Indeed, in the general case, it is not easy to obtain any performance guarantees for such algorithms. 

The next section puts forth a regularized version of \eqref{main} that is amenable to an online algorithms. A general online optimization framework is developed that can be utilized to construct online solutions to such tracking problems with arbitrary convex objectives, subject to certain regularity conditions. The framework will also allow us to develop bounds on the cumulative difference between the corresponding offline and online solutions. Different from the MPC-based Algo.~\ref{algo1}, the online algorithm does not suffer from infeasibility issues, and the resulting performance guarantees can be utilized towards pre-selecting the problem parameters.
\section{Online Algorithms}
This section develops online algorithms for solving the D2D trajectory optimization problem. As in Sec. \ref{offline}, it is assumed that the users are non-cooperative so that the current location of user~2 given by $\x_2(t)$ and the current destination of user~1 denoted by $\d_1(t)$ are exogenous and only available to user~1's device at time $t$. For the sake of simplicity, we will henceforth consider the utility function $U(x) = -x^2$. The extension to arbitrary concave utilities is left as future work, and is not pursued here. Since only the trajectory of user~1 is being optimized, we drop the subscripts from the variables and parameters corresponding to user~1. That is, the new variables and parameters of user~1 become $\x(t)$, $v$, $T$, $\delta$, $\s$ and $\d(t)$. The section is split into two parts: general results for arbitrary convex functions are developed first and application to the D2D problem at hand is discussed next. 
\subsection{General problem and the associated guarantees}
In this subsection, we consider the following general optimization problem 
\begin{align}\label{online}
\{\hat{\x}^r(t)\}_{t=2}^{T'} &= \arg\min_{\{\x(t)\}_{t=2}^{T'}} \sum_{t=1}^{T'} f_t(\x(t)) \\
& \text{s. t. } \norm{\x{(t+1)}-\x(t)}\leq v\nonumber\\
&\text{for all } \ 1\leq t < T'.\label{const_online}
\end{align}
where $T':=T+\delta$ and the initial position $\x(1) = \s$ is fixed. Different from \eqref{main}, the constraint to reach the final destination is relaxed but is incorporated within the objective function. Specifically, we utilize a regularization approach, and the modified objective function can take the form:
\begin{align}
f_t(\x(t))&:=\lambda(t)\left\| \x(t)-\x_2(t)\right\|^2+(1-\lambda(t))\left\| \x(t)-\d(t)\right\|^2 \nonumber\\
&= \norm{\x(t)-\bell(t)}_2^2 \label{loss}
\end{align}
where we define the leading path $\bell(t) := \lambda(t) \x_2(t) + (1-\lambda(t))\d(t)$. The sequence $\{\lambda(t)\}$ is included in order to control the relative importance assigned to being near user~2 and reaching the destination. Specifically, it is required that $\lambda(t)$ is a decreasing sequence that goes from 1 to 0 as $t$ goes from $1$ to $T'$. Such a choice of $\lambda(t)$ ensures that the user places increasingly higher importance to reaching the final destination. Note that since $\lambda(t) \in [0,1]$, the optimization problem in \eqref{online} is still convex. 

It is remarked that the solution to \eqref{main} and \eqref{online} will generally not be the same. In particular, since the constraint in \eqref{dest_rel} is relaxed, the trajectory in \eqref{online} might stop short of the destination, depending on the manner in which $\lambda(t)$ is decreased. It will however be shown in Sec. \ref{simulations} that $\norm{\x(T+\delta)-\d}$ is generally small for the choice $\lambda(t) = t/(T+\delta)$. We emphasize that the cost function need not be the least-squares loss function in \eqref{loss}. More generally, it may be reasonable to choose convex functions of the form $f_t(\x(t)) := g(\norm{\x-\bell(t)})$ for some scalar function $g$.  


Towards solving \eqref{online} in an online manner, we make use of the online gradient descent (OGD) algorithm which results in the following update 
\begin{align}
\hat{\x}^o(t+1)=& \hat{\x}^o(t)-\frac{1}{\gamma}\nabla f_t(\hat{\x}^o(t))\label{Algo}
\end{align}
where $\gamma$ is the learning rate {independent of $T'$} and must be chosen according to certain rules that we will discuss later. In present scenario, we also require that the problem parameters and the objective function be chosen so as to adhere to the following constraints. 
\begin{itemize}
	\item [(\textbf{A1})] \textbf{Strong convexity:} The time-varying objective function $f_t(\x)$ is strongly convex with parameter $\mu$, i.e., the function $f_t(\x) - \frac{\mu}{2}\norm{\x}^2_2$ is a convex function. 
	\item [(\textbf{A2})] \textbf{Lipschitz continuous gradient:} The gradient of the time varying objective function $f_t(\x)$ is bounded as $\norm{\nabla f_t(\x)}\leq G$ and Lipschitz continuous, i.e., $\norm{\nabla f_t(\x)-\nabla f_t(\y)}\leq L \norm{\x-\y}$ for all $t$. 
	\item [(\textbf{A3})] \textbf{Bounded variations: } The OGD iterates adhere to the velocity constraint $\norm{\hat{\x}^o(t+1)-\hat{\x}^o(t)}\leq v$ for all $t\geq 1$. 
\end{itemize}
Of these, the strong convexity and Lipschitz gradient constraints are standard and can be satisfied via appropriate choice of the cost function \cite{mokhtari2016online}. The bounded iterate variation assumption (\textbf{A3}) is however not standard and must be explicitly checked a priori. In general (\textbf{A3}) holds if the function gradient is bounded and $\gamma \geq G/v$. 

In order to quantify the performance of the algorithm, we will use tools from the online learning framework. Within this context, the user is treated as a learner that takes an action $\hat{\x}^o(t)$ at time $t$. In response to the user's action, the adversary selects the function $f_{t}(\cdot)$ and reveals $\nabla f_{t}(\hat{\x}^o(t))$ to the learner, and so on. The performance of the learner is often measured through the notion of regret. Specifically, we define the regret for our problem as
\begin{align}\label{regret}
\textbf{Reg}_{T'}:=\underbrace{\left[\sum_{t=1}^{T'} f_t(\hat{\x}^o(t))\right]}_{\text{online}}-\underbrace{\left[\sum_{t=1}^{T'} f_t(\hat{\x}^r(t))\right]}_{\text{offline}}
\end{align}
which is the cumulative difference between the objective function evaluated at the optimal offline trajectory and the objective function  evaluated at iterates given by the online algorithm in \eqref{online}. We remark that the definition of regret in \eqref{regret} is motivated from the offline regret introduced in \cite{chen2017online} and is stronger than the more commonly used dynamic regret. Indeed, the dynamic regret is commonly defined as \cite{mokhtari2016online} 
\begin{align}\label{regret2}
\textbf{Reg}^D_{T'}:=\sum_{t=1}^{T'} [f_t(\hat{\x}^o(t))- \min_{\x(t)} f_t(\x(t))]
\end{align}
and allows the adversary to choose $\x(t)$ as the minimizer of $f_t(\cdot)$. In contrast, the offline regret in \eqref{regret} allows the adversary to choose the full trajectory in an offline manner as in \eqref{online}. Indeed, the constraint in \eqref{const_online} couples the values of $\hat{\x}^r(t)$ over time, and as a result the offline regret is always higher than the dynamic regret. Clearly, both the definitions coincide if there are no temporally coupled constraints in \eqref{online}. 


The regret bounds will be calculated in terms of the squared path length of the adversary, defined as
\begin{align}
S_{T'}^\star:=\sum_{t=1}^{T'-1}\norm{\hat{\x}^r(t+1)-\hat{\x}^r(t)}^2
\end{align}
Having defined the squared path length, we are ready to state the main result of this subsection, whose proof is deferred to Appendix \ref{proof_thm_1}. 
\begin{thm}\label{theorem1}
	Under the assumptions \textbf{(A1)-(A3)}, for $\gamma \geq L$, and $T':=T+\delta$, the sequence of $\x^o(t)$ generated by the algorithm in \eqref{algo_main} adheres to the regret bound 
	\begin{align}
	\textbf{Reg}_{T'}\leq \sqrt{T'\mathcal{O}(S_{T'}^\star)}.
	\end{align}
\end{thm}
The result in Theorem 1 states that for large value of $T'$ and for a sublinearly time-varying adversary, the online algorithm incurs a sublinear regret over the offline solution. It is remarked that this result is stronger than similar results that have been shown to hold for dynamic regret in \cite{mokhtari2016online}. Likewise, the offline regret is also known to be linear in general \cite{chen2017online}. In contrast, Theorem \ref{theorem1} establishes a generic sublinear bound on the offline regret under relatively mild assumptions. 

Before concluding, it is remarked that while we skipped the constraint $\x \in \S$ for ease of exposition, Theorem \ref{theorem1} continues to hold for the projected OGD $\hat{\x}^o(t+1)= P_{\S}\left(\hat{\x}^o(t)-\frac{1}{\gamma}\nabla f_t(\hat{\x}^o(t))\right)$ where $P(\S)$ denotes the projection operation onto the convex set $\S$. The next section discusses the application of the developed results to the problem at hand. 

\subsection{D2D trajectory problem}
As detailed earlier, the OGD algorithm can be applied to the D2D trajectory optimization problem using any convex function of the form $g(\norm{\x(t)-\bell(t)})$. For instance, the square loss function that results in the gradient being proportional to $\x(t)-\bell(t)$. At the start, when the user is far from $\bell(t)$, the user may take larger steps, while it may slow down as it gets closer to $\bell(t)$. 

Another interesting choice is motivated from the Huber function, and takes the form
\begin{align}
g(d) = \begin{cases} \frac{1}{2}d^2 & d \leq v\\
v(1-\mu)d + \frac{\mu}{2}d^2 - \frac{(1-\mu^2)v^2}{2} & d > v
\end{cases}
\end{align}
for some parameter $\mu$. Observe that the Huber function is same as the squared distance as long as the distance is less than $v$. However, whenever the distance is larger than $v$, the penalty is a convex combination of the linear and squared error penalties. The constants are adjusted to ensure that the function $f_t = g(\norm{\x(t)-\bell(t)})$ adheres to (\textbf{A2}). As compared to the squared loss function, the Huber loss puts a smaller penalty when the user is far from $\bell(t)$. 

For this choice of objective function, its gradient can be written as
\begin{align}
\nabla f_t(\x)&= \begin{cases} \x-\bell(t) & \hspace{-1.5cm} \norm{\x-\bell(t)} \leq v\\
v(1-\mu)\frac{\x-\bell(t)}{\norm{\x-\bell(t)}} + \mu(\x-\bell(t))
\end{cases}\\
\nabla f_t(\hat{\x}^o(t)) &=\mu(\hat{\x}^o(t)\!-\!\bell(t)) \!+\! (1\!-\!\mu)P_v(\hat{\x}^o(t)\!-\!\bell(t))
\end{align}
where the projection operation is defined as
\begin{align}
P_v(\w)  &:= \arg\min_{\v} \norm{\v-\w}^2
&\text{s. t.} \norm{\v} &\leq v.
\end{align}
Therefore the OGD updates take the form
\begin{align}
&\hat{\x}^o{(t\!+\!1)} = \hat{\x}^o(t) \!-\! \frac{\mu}{\gamma}(\hat{\x}^o(t)\!-\!\bell(t)) \!-\! \frac{1\!-\!\mu}{\gamma} P_v(\hat{\x}^o(t)-\bell(t))\label{algo_main}
\end{align}
A special case occurs when $\mu \approx 0$ and $\gamma = 1$, for which case, the updates become 
\begin{align}\label{proposed}
\hat{\x}^o(t+1) = \hat{\x}^o(t) - P_v(\hat{\x}^o(t)-\bell(t)).
\end{align}

It can be seen that the function gradient is Lipschitz continuous with parameter $L = 1$ and strongly convex with parameter $\mu$. Moreover, if the user operates in a compact region of diameter $R$, the gradient is bounded as $\norm{\nabla f_t(\x)}_2 \leq \mu R + v(1-\mu)$. Consequently, (\textbf{A3}) is satisfied if $\mu R + v(1-\mu) \leq \gamma v$ or equivalently, if $\gamma \geq \frac{\mu (R-v)+v}{v} \geq 1+\frac{\mu R}{v}$. For instance when $\mu$ is close to zero, we simply require $\gamma \geq 1$. 

\section{Numerical results}\label{simulations}
This section provides detailed simulations for the various formulations provided here. Regardless of the algorithm used, the different formulations and settings will be compared on the basis of the cumulative D2D rates. At each time slot $t$, the maximum achievable rate given by 
\begin{align}
R(t) = W\log \left(1+\frac{P_s(t)}{P_s(t)+\sigma^2}\right)\label{rate}
\end{align}
where $P_s(t)=\norm{\x_1(t)-\x_2(t)}^{-\alpha}$ is the scaled path loss component and $\sigma^2$ is the appropriately scaled noise power. As remarked earlier the achievable rate does not represent the actual rate seen by the users, but is used here only for relative performance evaluation. 
\begin{figure}
	\centering
	\hspace{0mm}\includegraphics[scale=0.26]{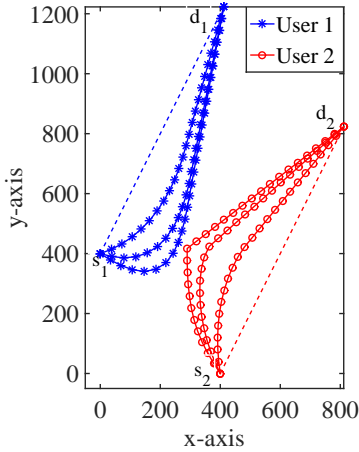}
	\captionsetup{font=scriptsize}
	\caption{D2D cooperative users offline approach for $\delta=1,3$ and $5$, trajectories becomes closer as $\delta$ increases and achieves the average data rate of $1.9$Mbps, $2.8$Mbps, and $3.5$ Mbps as compared to $1.1$Mbps for the direct path. }
	\label{first}
	\vspace{-5mm}
\end{figure}
We begin with the offline problem formulated in \eqref{main} for the case when $\S$ is simply a box in $\mathbb{R}^2$ and the utility function is $U(x) = -x^2$. We consider two pedestrians that start at the coordinates $\s_1 = [0~400]^T$ and $\s_2 = [400~0]^T$ and have destinations at $\d_1 = [400~1200]^T$ and $\d_2 = [800~800]^T$, respectively. Both users walk with the maximum speed of $v_1 = v_2 = 1$ m/s and take about $15$ minutes to reach their destinations through the direct path. In order to keep the problem sizes small, we divide this time into $T_1 = T_2 = 24$ time slots. As explained earlier, both the users are willing to incur an excess delay $\delta$ in reaching their destinations. Fig.~\ref{first} shows a plot of the user trajectories obtained from the offline algorithm in \eqref{main}. In the figure, the direct paths are shown as dashed lines while the optimal trajectories for different $\delta$ values are shown in color and with markers. Interestingly, the solution to the offline problem confirms to our intuition that the users first come close to each other and then part ways towards their respective destinations. Further intuition can be obtained by looking at the average D2D rate as a function of the allowable delay $\delta$, as shown in Fig. \ref{data_1}. The second y-axis of this figure depicts the increase in the total downloaded file size with respect to delay.  The other parameters used in Fig. 2 are $\alpha = 2.5$, $W = 10$ MHz, and $\sigma^2 = 0.2$. As expected, the average achievable rate continues to increase with $\delta$, as it allows larger deviations from the direct path. 
\begin{figure}
	\centering
	\includegraphics[width=0.89\linewidth, height = 0.5\linewidth]	{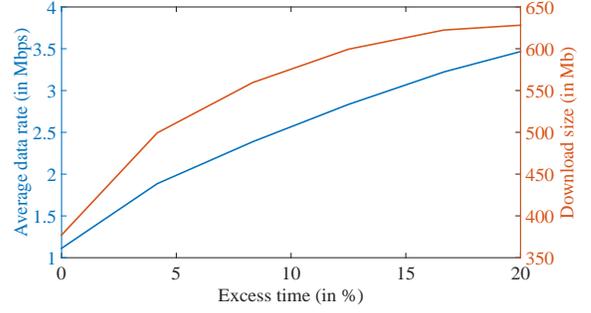}
	\captionsetup{font=scriptsize}
	\caption{Average data rate obtained with respect to increasing delay $\delta$ (excess time) for D2D cooperative users, and corresponding size of downloaded file. }
	\label{data_1}
\end{figure}

\begin{figure}
	\centering
	\hspace{0mm}\includegraphics[scale=0.3]{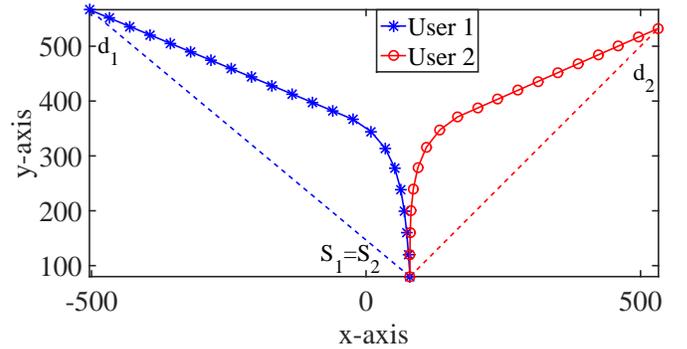}
	\captionsetup{font=scriptsize}
	\caption{D2D cooperative users with same starting and opposite destinations for $\delta=2$.  }
	\label{second}
	\vspace{0mm}
\end{figure}

In order to further examine the nature of the optimal trajectory, consider a special case when both users have the same starting point $\s_1 = \s_2 = [80\  80]^T$, but different destinations $\d_1 = [-400\ 480]^T$ and $\d_2 = [600 \ 600]^T $. The optimal trajectories are shown in Fig. \ref{second} for $\delta = 2$. In this case, the trajectory exhibits a knee region, where the users initially stay close together and then abruptly split up to head towards their respective destinations. 

\begin{figure}
	\centering
	\hspace{0mm}\includegraphics[scale=0.3]{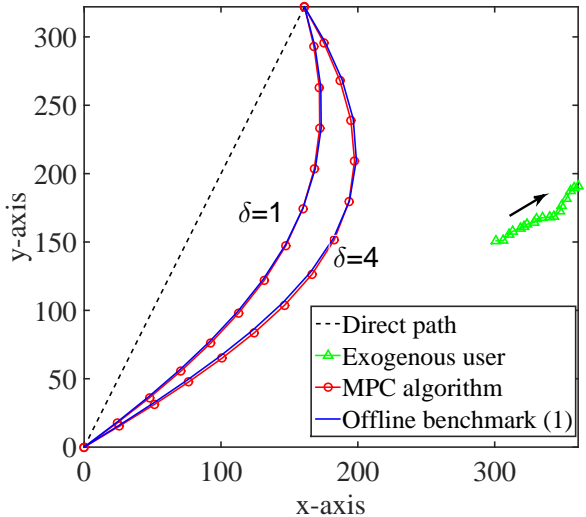}
	\captionsetup{font=scriptsize}
	\caption{D2D non cooperative users trajectory comparison for Offline methods, average data rate for direct path is $3.1$ Mbps, for MPC algorithm data rate is $3.357$ Mbps ($\delta=1$), $4.072$ Mbps ($\delta=4$), and for offline optimal data rate is  $3.359$ Mbps ($\delta=1$), $4.074$ Mbps ($\delta=4$) }
	\label{third}
\end{figure}

Next, we consider the non-cooperative case, where the users reveal their current locations but keep their future locations private. For all these cases, we assume the destination of user~1 to be fixed but allow the trajectory of user~2 to be exogenous.  As before we consider the square-law cost function and maximum velocity of $v = 1$ m/s. The user starts at the origin and heads towards the destination $\d = [150\ 300]^T$, taking about $6$ minutes to reach its destination. Time is discretized into slots such that $T =24 $ and the trajectory is shown for $\delta = 1$ and $\delta = 4$.  The offline benchmark scheme in Fig.~\ref{third} refers to the optimization problem in \eqref{main} with $\x_2(t)$ treated as an exogenous variable that is known in advance (arrow in Fig.~\ref{third} describes the direction of movement). The MPC-like method in Algo.~\ref{algo1} is plotted to account for the more realistic scenarios where $\x_2(t)$ is not known in advance and is only revealed in real-time. As evident from Fig. \ref{third}, the trajectory of Algo.~\ref{algo1} is quite close to that of the offline algorithm, since the speed of $\x_2(t)$ is much smaller than $v$. It was observed however that Algo.~\ref{algo1} is not very flexible, and suffers from infeasibility problems when the exogenous user moves at a higher speed.

\begin{figure}
	\centering
	\hspace{0mm}\includegraphics[scale=0.3]{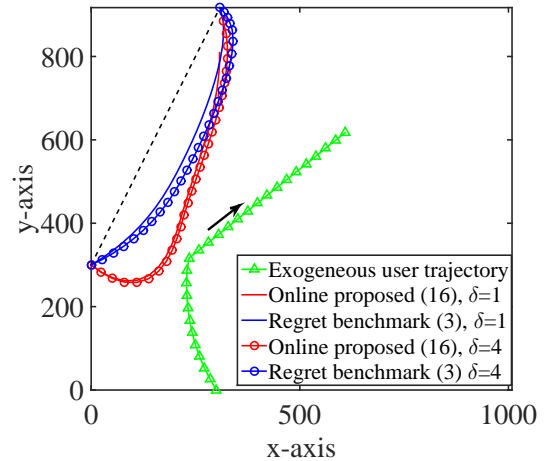}
	\captionsetup{font=scriptsize}
	\caption{Proposed online method compared with regret benchmark of \eqref{online} for $\delta=1$ and $\delta=4$ (marked lines).}
	\label{compa}
\end{figure}

\begin{figure}
	\centering
	\hspace{0mm}\includegraphics[width=\linewidth, height = 0.5\linewidth]{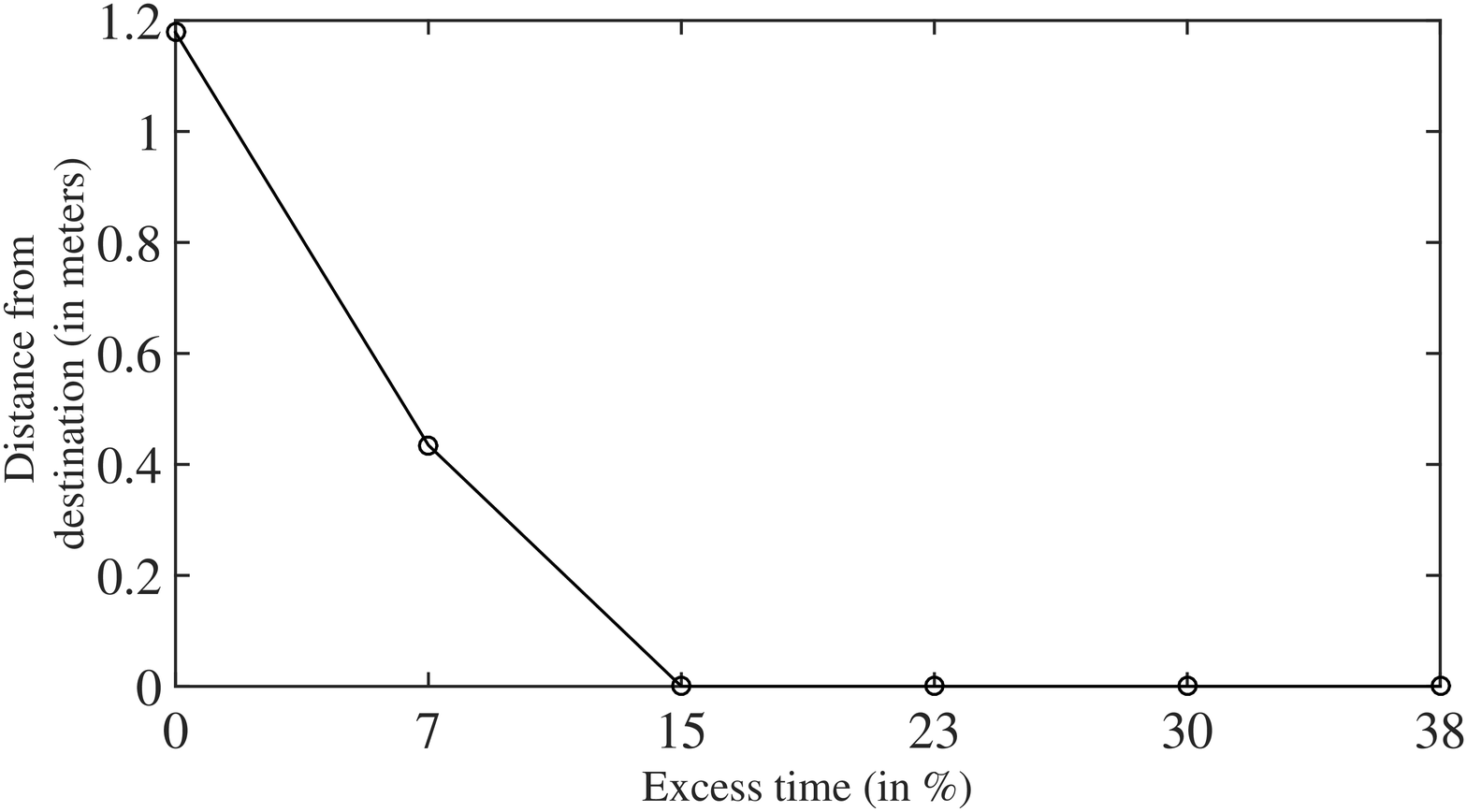}
	\captionsetup{font=scriptsize}
	\caption{Distance from the destination.}
	\label{distance_dest}
	\vspace{-3mm}
\end{figure}

\begin{figure}
	\centering
	\includegraphics[width=0.8\linewidth, height = 0.5\linewidth]
	{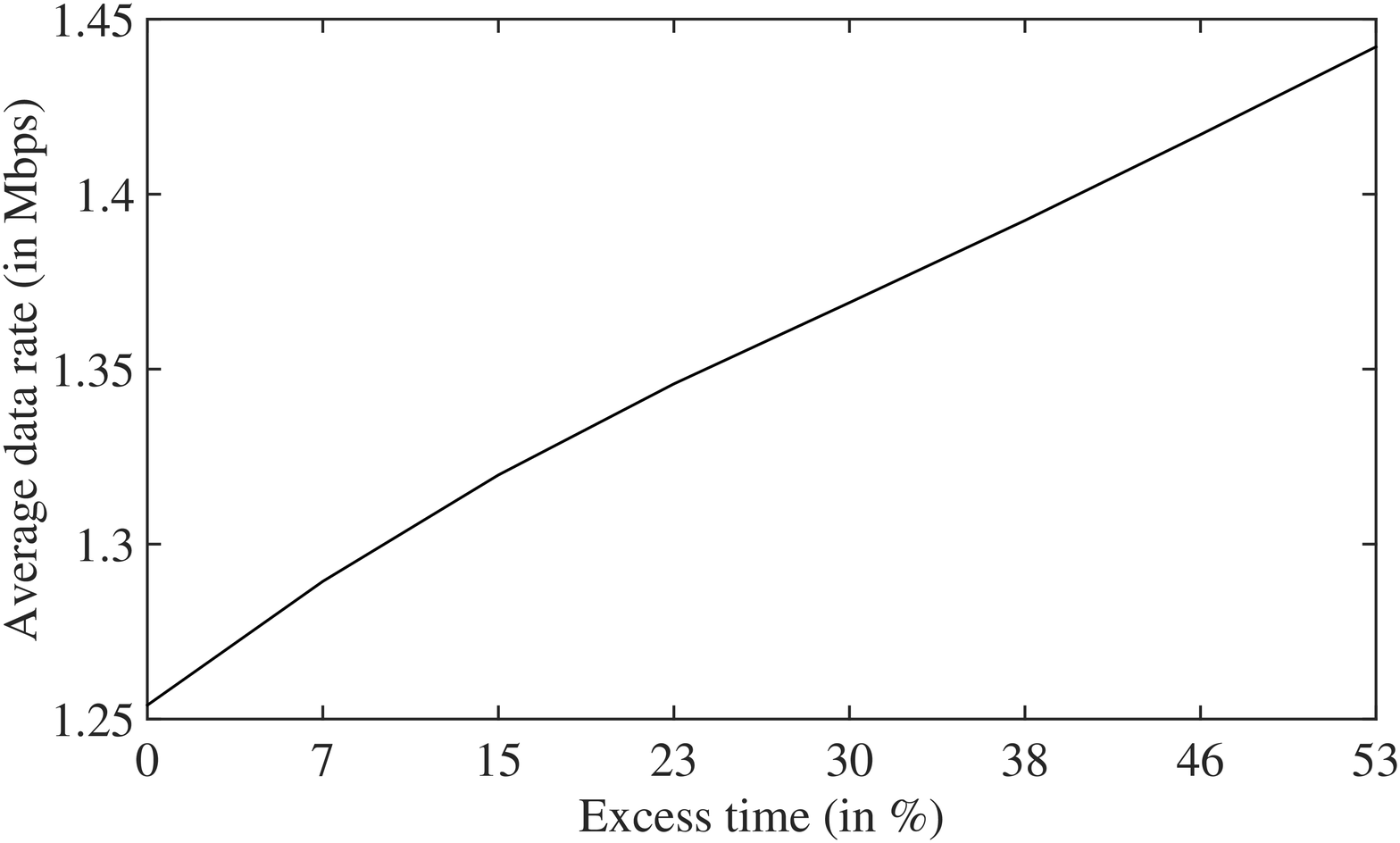}
	\captionsetup{font=scriptsize}
	\caption{D2D non-cooperative user performance in terms of average data rate achieved with respect to excess time by using proposed online algorithm of \eqref{proposed}.}
	\label{rate_2}
\end{figure}

Next, we present simulations for the online framework proposed in \eqref{proposed} for non-cooperative settings with $\mu=10^{-3}$ and $\gamma=1$. The online trajectory obtained for the proposed algorithm in  \eqref{proposed} is shown in Fig.~\ref{compa} for different values of delay $\delta$. The trajectory is compared with the offline regret benchmark defined in \eqref{online}. To plot this figure, the regularizer value used is $\lambda(t)=\frac{t}{T+\delta}$. If $\lambda(t)$ is not appropriately chosen, it may happen that the user will not be able to reach the destination. But this condition do not arises for sufficiently high $\delta$ as depicted in Fig.~\ref{distance_dest}. Note that for smaller values of $\delta$ (less excess time),  user is at a small distance from its final destination, which  reduces to zero for larger values of delays (more excess time) and user will always reach the destination. Further in Fig.~\ref{rate_2}, average data rate achieved by the proposed  online algorithm is shown against the excess time. 

\section{Conclusion and open problems}

This paper considered the problem of designing the trajectory of a pair of device-to-device users. The problem is formulated for a pair of users and is shown to be solvable in an offline manner. Motivated by the users need to keep their future locations private, we move on to develop online algorithms for the same. Specifically, we develop a modified OGD algorithm that incurs a sublinear offline regret, a result that has not been reported in the literature. Detailed simulations have been carried out to demonstrate the efficacy of the algorithm.

Since this is the first D2D trajectory optimization formulation of its kind, several aspects of the problem have been ignored or simplified, allowing us to build general-purpose and sophisticated approaches. In a realistic setting, it may be necessary to take care of constraints such as those arising from the presence of obstacles and roads. The present formulation cannot directly incorporate roads and obstacles since the bounds developed here rely on the area being convex. Next, the algorithms proposed here rely heavily on the availability of the exact location of the two users. In practice, if the privacy concerns prevent the users from sharing their location information, one could still attempt the trajectory optimization by sampling the received signal strength. However, the resulting uncertainty in the location information might require us to formulate the problem within a stochastic setting, and develop corresponding stochastic regret bounds. Another open problem that is fairly common is the case of three or more users, possibly interacting during different parts of the commute. For instance the user~1 may communicate with user~2 for the first 10 minutes and subsequently switch to a new user~3.

\appendices
\vspace{-3mm}
\section{Proof of Theorem~\ref{theorem1}}\label{proof_thm_1}
We begin with the following quadratic lower bound that can be obtained as a consequence of the strong convexity assumption in (\textbf{A1}) 
\begin{align}
f_t(\hat{\x}^r(t)) \geq& f_t(\hat{\x}^o(t)) + \nabla f_t(\hat{\x}^o(t))^T(\hat{\x}^r(t)-\hat{\x}^o(t)) \nonumber\\
&+ \frac{\mu}{2}\norm{\hat{\x}^r(t)-\hat{\x}^o(t)}_2^2
\end{align}	
Adding and subtracting $\nabla f_t(\hat{\x}^o(t))^T\hat{\x}^o(t+1)$ on the right hand side, we obtain
\begin{align}\label{strong}
f_t(\hat{\x}^r(t))- &\frac{\mu}{2}\norm{\hat{\x}^r(t)-\hat{\x}^o(t)}_2^2\\
&\geq f_t(\hat{\x}^o(t)) + \nabla f_t(\hat{\x}^o(t))^T(\hat{\x}^o(t+1)-\hat{\x}^o(t)) \nonumber\\
& +\nabla f_t(\hat{\x}^o(t))^T(\hat{\x}^r(t)-\hat{\x}^o(t+1)) \nonumber \\
&= f_t(\hat{\x}^o(t)) + \nabla f_t(\hat{\x}^o(t))^T(\hat{\x}^o(t+1)-\hat{\x}^o(t)) \nonumber\\
& +\gamma (\hat{\x}^o(t+1)-\hat{\x}^o(t))^T(\hat{\x}^o(t+1)-\hat{\x}^r(t)) \nonumber
\end{align}	
where the last equality follows from the the update in \eqref{proposed} for the general case. 
Note that the last term can be written as
\begin{align}
&\gamma (\hat{\x}^o(t+1)-\hat{\x}^o(t))^T(\hat{\x}^o(t+1)-\hat{\x}^r(t)) \\
&= \!\gamma (\hat{\x}^o(t\!\!+\!\!1)\!-\!\hat{\x}^o(t))^T(\hat{\x}^o(t)\!-\!\hat{\x}^r(t)) \!\!+\!\! \gamma \norm{\hat{\x}^o(t\!\!+\!\!1)\!-\!\hat{\x}^o(t)}_2^2\nonumber
\end{align}
The Lipschitz gradient condition yields
\begin{align}\label{Lips}
f_t(\hat{\x}^o(t+1)) \leq& f_t(\hat{\x}^o(t)) + \nabla f_t(\hat{\x}^o(t))^T(\hat{\x}^o(t+1)-\hat{\x}^o(t)) \nonumber\\
&+ \frac{L}{2}\norm{\hat{\x}^o(t+1)-\hat{\x}^o(t)}_2^2
\end{align}
Combining the two equations in \eqref{strong} and \eqref{Lips}, we obtain:
\begin{align}
f_t&(\hat{\x}^r(t))-f_t(\hat{\x}^o(t+1)) \nonumber\\
\geq& \frac{\mu}{2}\norm{\hat{\x}^r(t)-\hat{\x}^o(t)}_2^2  + \left(\gamma-\frac{L}{2}\right)\norm{\hat{\x}^o(t+1)-\hat{\x}^o(t)}_2^2 \nonumber\\
& + \gamma(\hat{\x}^o(t+1)-\hat{\x}^o(t))^T(\hat{\x}^o(t)-\hat{\x}^r(t))
\end{align}
Taking sum over $t = 1, \ldots, T'$,  using the definition of $\{\hat{\x}^r(t)\}_{t=1}^{T'}$ which are optimal points, and after rearranging, we have that
\begin{align}\label{lower_bound}
\sum_{t=1}^{T'}& (\hat{\x}^o(t+1)-\hat{\x}^o(t))^T(\hat{\x}^r(t)-\hat{\x}^o(t))\\
\geq&\frac{\mu}{2\gamma}\sum_{t=1}^{T'}\norm{\hat{\x}^r(t)-\hat{\x}^o(t)}_2^2 \nonumber\\
&
+ \frac{1}{\gamma}\left(\gamma-\frac{L}{2}\right)\sum_{t=1}^{T'}\norm{\hat{\x}^o(t+1)-\hat{\x}^o(t)}_2^2. \nonumber
\end{align}
Further, for the first term, we have that 
\begin{align}
\norm{\hat{\x}^o(t\!\!+\!\!1)\!-\!\hat{\x}^r(t)}_2^2\!\!=& \!\norm{\hat{\x}^o(t\!\!+\!\!1)\!-\!\hat{\x}^o(t)}^2_2+\norm{\hat{\x}^r(t)\!\!-\!\!\hat{\x}^o(t)}^2_2\!\nonumber\\
&-2(\hat{\x}^o(t\!\!+\!\!1)\!-\!\hat{\x}^o(t))^T(\hat{\x}^r(t)\!-\!\hat{\x}^o(t)).
\end{align}
Take the summation from $t=1$ to $T'$ on both sides, we get
\begin{align}\label{above}
\sum_{t=1}^{T'}&\norm{\hat{\x}^o(t\!\!+\!\!1)\!-\!\hat{\x}^r(t)}_2^2\!\!\nonumber\\&= \sum_{t=1}^{T'}\!\norm{\hat{\x}^o(t\!\!+\!\!1)-\hat{\x}^o(t)\!}^2_2+\sum_{t=1}^{T'}\norm{\hat{\x}^r(t)-\hat{\x}^o(t)\!\!}^2_2\!\nonumber\\
&-2\sum_{t=1}^{T'}(\hat{\x}^o(t\!\!+\!\!1)-\hat{\x}^o(t)\!)^T(\hat{\x}^r(t)-\hat{\x}^o(t)\!).
\end{align}
Utilizing the lower bound in \eqref{lower_bound} for the last term on right hand side of \eqref{above}, we obtain
\begin{align}\label{upper_bound}
\sum_{t=1}^{T'} \norm{\hat{\x}^o(t\!\!+\!\!1)\!-\!\hat{\x}^r(t)}^2_2
&\!\leq\! \frac{\gamma\!-\!\mu}{\gamma}\!\sum_{t=1}^{T'}\norm{\hat{\x}^r(t)-\hat{\x}^o(t)}^2_2 \\
& - \frac{\gamma-L}{\gamma}\sum_{t=1}^{T'}\norm{\hat{\x}^o(t+1)-\hat{\x}^o(t)}^2_2.\nonumber
\end{align}
Next, it follows from  the use of  Peter-Paul inequality with parameter $\eta > 0$:
\begin{align}\label{peter}
&\norm{\hat{\x}^o(t+1)-\hat{\x}^r(t+1)}_2^2 \\
&\leq \!\!(\!1\!+\!\eta\!)\!\norm{\hat{\x}^o(t\!\!+\!\!1)-\hat{\x}^r(t)}_2^2 \!+\! \left(1\!\!+\!\!\frac{1}{\eta}\right)\!\norm{\hat{\x}^r(t+1)-\hat{\x}^r(t)}_2^2. \nonumber
\end{align}
Taking the summation over $t$ on both sides of \eqref{peter}, and then utilize the upper bound of \eqref{upper_bound}, we get 
\begin{align}
\sum_{t=1}^{T'} &\norm{\hat{\x}^o(t+1)-\hat{\x}^r(t+1)}_2^2 \nonumber\\
\leq &\frac{(1+\eta)(\gamma-\mu)}{\gamma}\sum_{t=1}^{T'}\norm{\hat{\x}^r(t)-\hat{\x}^o(t)}^2_2 \nonumber\\
& - \frac{(1+\eta)(\gamma-L)}{\gamma}\sum_{t=1}^{T'}\norm{\hat{\x}^o(t+1)-\hat{\x}^o(t)}^2_2 \nonumber\\
& + \left(1+\frac{1}{\eta}\right)\sum_{t=1}^{T'}\norm{\hat{\x}^r(t+1)-\hat{\x}^r(t)}_2^2.
\end{align}
Assuming that $\hat{\x}^o(1) = \hat{\x}^r(1)$ and $\hat{\x}^r(T'+1) = \hat{\x}^r(T')$, we obtain for $(1+\eta)(1-\frac{\mu}{\gamma}) \in (0,1)$:
\begin{align}
\sum_{t=1}^{T'}&\norm{\hat{\x}^o(t)-\hat{\x}^r(t)}^2_2 \nonumber
\\ & \leq \frac{(1+\eta)\gamma }{\eta(\mu(1+\eta)-\eta\gamma)}S_{T'}^\star - \frac{(1+\eta)(\gamma-L) }{\mu(1+\eta)-\eta\gamma}O_{T'}\nonumber\\
&-\norm{\hat{\x}^o(T'+1)-\hat{\x}^r(T'+1)}_2^2
\end{align}
where $O_{T'}:=\sum_{t=1}^{T'-1}\norm{\hat{\x}^o(t+1)-\hat{\x}^o(t)}^2_2$. The second term and the third term in the above expression can be dropped since $\gamma\geq L$. Further, minimizing the term $\frac{(1+\eta)\gamma }{\eta(\mu(1+\eta)-\eta\gamma)}$ associated with first term with respect to $\eta$, we obtain 
\begin{align}
\eta=\sqrt{\frac{\gamma}{\gamma-\mu}}-1
\end{align}
resulting in the following bound
\begin{align}\label{final}
\sum_{t=1}^{T'}\norm{\hat{\x}^o(t)-\hat{\x}^r(t)}^2_2 \leq \frac{S_{T'}^\star}{1-\sqrt{1-\mu/\gamma}}\approx\sqrt{\frac{2\gamma}{\mu}}S_{T'}^\star.
\end{align}
Since $\mu$ and $\gamma$ are constants not dependent on $T'$, the bound can be compactly written as $\sum_{t=1}^{T'}\norm{\hat{\x}^o(t)-\hat{\x}^r(t)}^2_2 \leq  \mathcal{O}(S_{T'}^\star)$. Finally, we obtain the required regret bound by using first order convexity condition, Cauchy-Schwartz inequality, and gradient boundedness as follows 
\begin{align}\label{bound}
\sum_{t=1}^{T'}f_t(\hat{\x}^o(t))-f_t(\hat{\x}^r(t))&\leq \sum_{t=1}^{T'}\nabla f_t(\hat{\x}^o(t))^T(\hat{\x}^o(t)-\hat{\x}^r(t))\nonumber\\
&\leq G\sum_{t=1}^{T'}\norm{\hat{\x}^o(t)-\hat{\x}^r(t)}_2.
\end{align}
It holds that 
\begin{align}\label{square_ine}
\sum_{t=1}^{T'}\norm{\hat{\x}^o(t)-\hat{\x}^r(t)}_2\leq \sqrt{T'\sum_{t=1}^{T'}\norm{\hat{\x}^o(t)-\hat{\x}^r(t)}^2_2 }.
\end{align}
Utilizing \eqref{square_ine} and  \eqref{final} into \eqref{bound}, we get
\begin{align}
\sum_{t=1}^{T'}f_t(\hat{\x}^o(t))-f_t(\hat{\x}^r(t))\leq \sqrt{T'\mathcal{O}(S_{T'}^\star)}.
\end{align}

Hence proved.

										\footnotesize
										\bibliographystyle{IEEEtran} 
										\bibliography{IEEEabrv,references}
										
									\end{document}